\documentclass[11pt]{article}
\usepackage{stmaryrd}
\usepackage{amssymb}
\usepackage{color, amsmath,amssymb, amsfonts, amstext,amsthm, latexsym}
\usepackage{amssymb, epsfig, amssymb, latexsym}
\usepackage{amsmath}
\usepackage{graphicx}
\usepackage{longtable}
\usepackage{cite}
\usepackage{subfigure}
\numberwithin{equation}{section} \topmargin -0.4in

\renewcommand{\phi}{\varphi}

\theoremstyle{definition}

\theoremstyle{remark}

\numberwithin{equation}{subsection}

\title{Post-grazing dynamics of a vibro-impacting energy generator}

\author{Larissa Serdukova\thanks{School of Mathematics, Georgia Institute of Technology, Atlanta, GA 30313, USA. E-mail: larissa.serdukova@math.gatech.edu}, Rachel Kuske\thanks{School of Mathematics, Georgia Institute of Technology, Atlanta, GA 30313, USA. E-mail: rachel@math.gatech.edu}, Daniil Yurchenko \thanks{IMPEE, Heriot-Watt University, Edinburgh, EH14 4AS, UK. E-mail: d.yurchenko@hw.ac.uk}.}   

\begin{document}
\maketitle
\begin{abstract}
The motion of a forced vibro-impacting inclined energy harvester is investigated in parameter regimes with asymmetry in the number of impacts on the bottom and top of the device. This motion occurs beyond a grazing bifurcation, at which alternating top and bottom impacts are supplemented by a zero velocity impact with the bottom of the device. For periodic forcing, we obtain semi-analytical expressions for the asymmetric periodic motion with a ratio of 2:1 for the impacts on the device bottom and top, respectively. These expressions are derived via a set of  nonlinear maps between different pairs of impacts, combined with impact conditions that provide jump discontinuities in the velocity. Bifurcation diagrams for the analytical solutions are complemented by a linear stability analysis around the 2:1 asymmetric periodic solutions, and are validated numerically. For smaller incline angles, a second grazing bifurcation is numerically detected, leading to a 3:1 asymmetry. For larger incline angles, period doubling bifurcations precede this bifurcation. The converted electrical energy per impact is reduced for the asymmetric motions, and therefore less desirable under this metric.  
\end{abstract}



\section{Introduction}

Energy Harvesting (EH) from ambient vibrations was proposed almost two decades ago as an attractive alternative to power supplies or as renewable sources of energy for rechargeable batteries. Since then the gaps in the linear theory of EH have been filled  with different methods of energy conversion, based on single-degree-of freedom, multi-degree-of freedom and/or continuous (rods and beams) linear systems on the nano, micro and macro scales \cite{Stephen2006,Bowen2016,Briand2015,Yang2014, Elvin2013}. The excitement regarding the potential of linear EH systems has significantly decreased since then due to low energy densities of the linear devices, narrow bandwidth and high natural frequency in nano- and micro-scale systems, which are difficult to realize in many practical applications. These and other adverse factors lead to insufficient power output to power or recharge a battery. The deficiencies in the development of linear EH devices has slowed the proliferation of wireless sensors, particularly critical in the Internet of Things paradigm.   

The above limitations in the linear theory of EH have motivated wide-spread efforts on parametrically excited, nonlinear and non-smooth systems. The idea behind parametrically excited systems is the use of large system responses near instabilities, e.g. see \cite{Wiercigroch2011,Kecik2013,Jia2014,Bobryk2016,Bobryk2016a,Yerrapragada2017,Yurchenko2018,Dotti2019,Kuang2019,Avanco2019}, among others. Within the huge range of nonlinear EH systems, there are some particular themes of note; natural single-potential nonlinearities (classical continuous nonlinear systems like the Van-der-Pol or Duffing oscillator, a pendulum, etc. as in \cite{Ghouli2017,Lingala,Zhu2015,Quinn2011,Sebald2011,Green2012,Daqaq2014}), natural or imposed geometrical nonlinearities (systems with double, triple or multiple stable equilibriums, as in \cite{Pellegrini2013,Gao2016,Fu2019,Hawes2017,Zhou2018,Huang2019,Deng2019,Yang2019}), systems with a nonlinear interaction such as flow-induced vibration systems (see \cite{Abdelkefi2012,Wang2019,Wang2019a,Jin2019,Orrego2017} and references therein), and systems with strongly nonlinear or discontinuous nonlinearities like dry friction, piecewise discontinuity or vibroimpacts \cite{Ibrahim2009,Dimentberg2004}. It has been shown that the nonlinear mechanisms for EH are far more beneficial than linear ones. This observation follows from the typical structure of the response amplitude vs. forcing frequency or backbone curve, showing a wider bandwidth with higher response amplitude away from a main resonance frequency. However, the design and optimization of a nonlinear energy harvester is far more complex, with limited explicit analytical results, thus requiring extensive complementary experiments or numerics. The available approximation techniques can estimate the response within only a narrow range of parameters imposed by the mathematical assumptions necessary for the applied averaging procedure, typically based on a weakly nonlinear model with  small forcing. 

Vibro-impact systems have rich phenomenological behaviors, manifesting various nonlinear phenomena like bifurcations, grazing and chaos. These effects have been studied in deterministic and stochastic vibro-impact systems, as in \cite{Nordmark1997,Chillingworth2002,Wagg2004,Wagg2005,Simpson2013,Kumar2016,DiBernardo2008,Luo2013,Luo2004,Simpson2018} among others. The models of vibro-impact systems include piecewise linear stiffness \cite{Shaw1983,ElAroudi2014} as well as rigid barriers and instantaneous impacts leading to a velocity jump for inelastic impacts. EH devices that utilize vibro-impact dynamics as a main energy absorption mechanism were developed and studied in a number of publications \cite{Borowiec2014,Gendelman2015,Truong2016,Bendame2016}. While often such systems are limited to computational results only, certain settings  allow an analytical or semi-analytical treatment when the motion is composed of a sequence of trajectories described (semi-)analytically. Such an approach translates the piecewise continuous behavior into a sequence of maps, amenable to analytical treatment \cite{Luo2013}. This methodology has certain benefits since it allows bifurcation and stability analyses of various periodic regimes that may occur in the system. Of course, for more complex motions a series of maps is necessary, making these derivations more tedious and cumbersome. 

Recently, \cite{Yurchenko2017} proposed a novel vibro-impact energy harvesting (VI-EH) device utilizing dielectric elastomeric (DE) membranes. There it was shown that the performance of such VI-EH depends strongly on the relationship between the excitation and device parameters, leading to various vibro-impact regimes with a low or high power output. The device consists of a forced cylinder  with a ball moving freely inside of it, impacting DE membranes covering both ends of the cylinder. Each membrane is composed of the DE material sandwiched between two compliant electrodes, acting as a variable capacitance capacitor. Impact of the DE membrane by the ball influences its motion while deforming the membrane, leading to energy harvesting via the properties of variable capacitance. We characterize the motion of the VI-EH in terms of the ratio n:m of impacts per period of the forcing, where $n$ and $m$ are the number of impacts against the bottom and top membrane, respectively. Here we restrict our attention to $m=1$, as these types of solutions appear over significant parameter ranges when the cylinder is inclined.

The bifurcations and linear stability of the 1:1 periodic motion (two alternating bottom and top impacts per period) was investigated in \cite{Serdukova2019}, demonstrating the  influence of parameters such as the length of the cylinder, excitation parameters, and incline angle on this motion and the corresponding VI-EH power output. However, this study did not consider adjacent parameter regimes where n:1 periodic motions, period doubling bifurcation, and chaotic motion were observed numerically. In this paper we take a wider view to study  semi-analytical solutions and stability conditions for 2:1 (three impacts per period) periodic motion, which has implications for the device's energy harvesting potential. A mechanical model and equations of motion of the VI-EH are described in Section 2, together with a review of the results from \cite{Serdukova2019}. In Section 3 semi-analytical expressions for a 2:1 periodic motion are derived through three nonlinear maps, corresponding to the three impacts per period, combined with the impact conditions. A linear stability analysis for this motion is given in Section 4. The voltage output of the 2:1 periodic motion is shown in Section 5 and compared with that for the 1:1 periodic motion, together with comparisons of different metrics for the average energy available for harvesting. Finally, conclusions are drawn together with recommendations for the device design.                

\section{Previous results for the vibro-impacting energy harvester (VI-EH)}\label{sec:previous}

\begin{figure}[ht]
\centering
\scalebox{0.45}
{\includegraphics{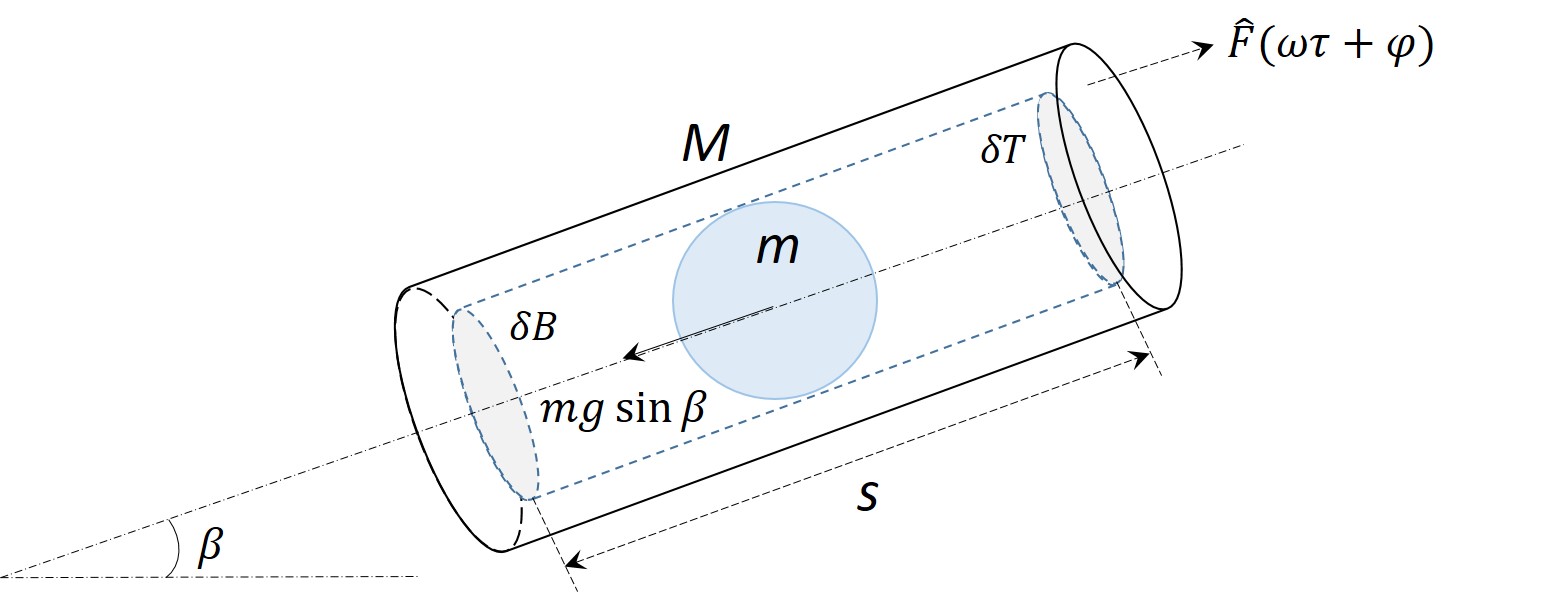}}
\caption{\scriptsize{A mechanical model of a vibro-impact energy harvester adapted from \cite{Yurchenko2017}.}}
\label{Fig00}
\end{figure}

We give the equations of motion for the inclined VI-EH with two DE membranes for harvesting energy from ambient vibrations as developed in \cite{Yurchenko2017}, as shown in the schematic of Figure \ref{Fig00}. 
The cylinder of mass $M$ and length $s$ is subject to a harmonic excitation $\hat{F}(\omega \tau + \varphi)$ with period $2 \pi/\omega$. Then the position of its center $X(\tau)$ satisfies
\begin{eqnarray}
\frac{\partial ^2{X}}{\partial \tau^2}=\frac{\hat{F}(\omega \tau + \varphi)}{M}. \label{cyl_dim}
\end{eqnarray} 
Between impacts, the  ball of mass $m$ ($M\gg m$) rolls freely inside of the cylinder driven only by  gravity ($g = 9.8 \; {\rm m/s}^2$  is the gravitational constant), with position $x$ given by
\begin{eqnarray}
\frac{\partial ^2{x}}{\partial \tau^2}
=-G = -g\sin\beta,   \label{ball_dim}
\end{eqnarray}
until it collides with one of the membranes causing its deformation. At impact, the velocity of the ball changes in sign and magnitude according to
\begin{eqnarray}\label{impact_dim}
\left(\frac{\partial  {x}}{\partial \tau}\right)^+=-r \left(\frac{\partial  {x}}{\partial \tau}\right)^- + (r+1) \frac{\partial  {X}}{\partial \tau}.
\end{eqnarray}
Here $r$ is the coefficient of restitution, and superscripts ${\ }^-$ and ${\ }^+$ indicate the velocities of the ball just before and after each impact, respectively. We assume that the velocity of the cylinder $\dot{X}$ does not change with an impact for $m$ negligible relative to $M$.

To track the dependence of periodic motions in terms of the parameters, it is valuable to use dimensionless equations of motion in terms of the relative variables. For this purpose we non-dimensionalize the original system \eqref{cyl_dim}, \eqref{ball_dim} with the substitutions
\begin{eqnarray}\label{nondXYt}
X(\tau)=\frac{\parallel \hat{F} \parallel \pi^2}{M \omega^2} \cdot X^*(t), \ \  \frac{\partial  {X}}{\partial \tau}=\frac{\parallel \hat{F} \parallel\pi}{M \omega} \cdot \dot{X}^*(t), \ \  \tau=\frac{\pi}{\omega} \cdot t \, ,
\end{eqnarray}
where  $\parallel \hat{F} \parallel$ is an appropriately defined norm of the strength of the forcing $\hat{F}$ and ''$\; \dot{\ }$ `` indicates the derivative with respect to $t$. Then the dimensionless equations of motion in terms of the relative position $Z(t)$ and velocity $\dot{Z}(t)$ are
\begin{eqnarray}
& & {Z}={X}^*-{x}^*, \qquad \dot{Z}=\dot{X}^*-\dot{x}^*\nonumber\\
& &\ddot{Z}=\ddot{X}^*-\ddot{x}^*= F(\pi t+ \varphi)+ \dfrac{M g \sin \beta}{\parallel \hat{F} \parallel}=f(t)+\bar{g}, \label{eq_relative}
\end{eqnarray}
where the non-dimensional forcing $F$ has the unit norm, i.e. $\parallel F \parallel = 1$, and period $2$. Then the impact condition \eqref{impact_dim} in terms of the non-dimensional relative variables for the $j$-th impact at time $t=t_j$ is
\begin{align}
&Z_j=X^*(t_j)-x^*(t_j)=\pm \frac{d}{2}, \ \
\mbox{ for } x \in \partial B \; (\partial T) \mbox{ the sign is } +(-) \nonumber \\
&\dot{Z}_j^+=-r \dot{Z}_j^-\, , \label{eq_impact}
\end{align}
where $d=\frac{s M \omega^2}{\parallel \hat{F} \parallel \pi^2}$ is the length of the cylinder, $\partial B$ and $\partial T$ are the bottom and top membranes of the energy harvesting system.  

By integrating \eqref{eq_relative} for $t\in (t_j, t_{j+1})$ and applying \eqref{eq_impact}, we obtain the expressions for the relative velocity and displacement between two impacts
\begin{align}
&\dot{Z}(t)=-r \dot{Z}^-_{j} + \bar{g}(t-t_{j}) + F_1(t)-F_1(t_{j}), \label{zdot_t}\\
&Z(t)=Z^-_{j}-r \dot{Z}^-_{j}(t-t_{j}) +\frac{\bar{g}}{2}(t-t_{j})^2+F_2(t)-F_2(t_{j})-F_1(t_{j})(t-t_{j}), \label{z_t}
\end{align}
where $F_1(t)=\int f(t)dt$ and $F_2(t)=\int F_1(t)dt$. In the following expressions, the superscripts ''${\ }^- \;$`` are omitted, since \eqref{zdot_t}-\eqref{z_t} are in terms $Z^-$ and $\dot{Z}^-$ only. Evaluating \eqref{zdot_t}-\eqref{z_t} at impact times $t=t_{j+1}$, we obtain equations defining the four basic nonlinear maps $P_l, \; l=1, 2, 3, 4$ for the corresponding transitions between impacts,
\begin{align}
&P_{1} \; : \; (Z_j\in \partial B, \dot{Z}_j, t_j) \; \mapsto \; (Z_{j+1}\in \partial B, \dot{Z}_{j+1}, t_{j+1}),\nonumber\\
&P_{2} \; : \; (Z_j\in \partial B, \dot{Z}_j, t_j) \; \mapsto \; (Z_{j+1}\in \partial T, \dot{Z}_{j+1}, t_{j+1}),\nonumber\\
&P_{3} \; : \; (Z_j\in \partial T, \dot{Z}_j, t_j) \; \mapsto \; (Z_{j+1}\in \partial B, \dot{Z}_{j+1}, t_{j+1}),
\end{align}
and similarly, for $P_4$ for the $\partial T \mapsto \partial T$ transition. Here we restrict our attention to $P_1$, $P_2$ and $P_3$ transitions, since only these play a role in the attracting 2:1 motion. The mathematical expressions for these maps take different forms depending on whether $Z_j$ and $Z_{j+1}$ are located on either $\partial B$ or $\partial T$. Specifically, for $t=t_{j+1}$, 
  \eqref{zdot_t} - \eqref{z_t} are given by
\begin{align}
&\dot{Z}_{j+1}=-r \dot{Z}_j + \bar{g} (t_{j+1}-t_j) + F_1(t_{j+1})-F_1(t_j),\label{zdotj_p1}\\
&D_\ell=-r \dot{Z}_j (t_{j+1}-t_j) +\frac{\bar{g}}{2} (t_{j+1}-t_j)^2 +F_2(t_{j+1})-F_2(t_j)-F_1(t_j) (t_{j+1}-t_j).\label{zj_p1}
\end{align}
where $D_1=D_4 =0$, $D_2 = -d$ and $D_3 = d$.

In \cite{Serdukova2019}, the expressions \eqref{zdotj_p1}-\eqref{zj_p1} for the maps $P_2$ and $P_3$  over the time intervals $(t_{k-1}, t_k)$ and $(t_k, t_{k+1})$  are combined with periodic and impact conditions to derive equations for the triples  $(\dot{Z}_{k},\varphi_{k},\Delta t_{k})$ corresponding to 1:1 periodic solutions. Throughout this paper $\Delta t_{k}=t_k-t_{k-1}$ for any $k$, and $\varphi_{k} =$mod$(\pi t_{k} + \varphi, 2\pi)$ is the phase shift of the $k^{\rm th}$ impact relative to that of the forcing $f(t)$. The resulting expressions provide the dependence of 1:1 motions on the combinations of the parameters $d$, $r$, $\bar{g}$ and $\hat F$. A calculation of the energy output follows directly from these triples, based on the deformation of the membrane that depends explicitly on $\dot{Z}_k$. Given a constant input voltage, $U_{\rm in}$, applied to the membranes,  the change in charge $Q$ across the capacitor is given by $\Delta Q=U\Delta C$, where $C$ is its capacitance. The charge $Q$ increases as the membrane's shape is  restored,  producing a higher voltage $U_k$ at the $k^{\rm th}$ impact, with resulting energy to be harvested in terms of the difference, $U_k- U_{\rm in}$, which we refer to as the output voltage.

\begin{figure}[t!]
\centering
\scalebox{0.60}
{\includegraphics{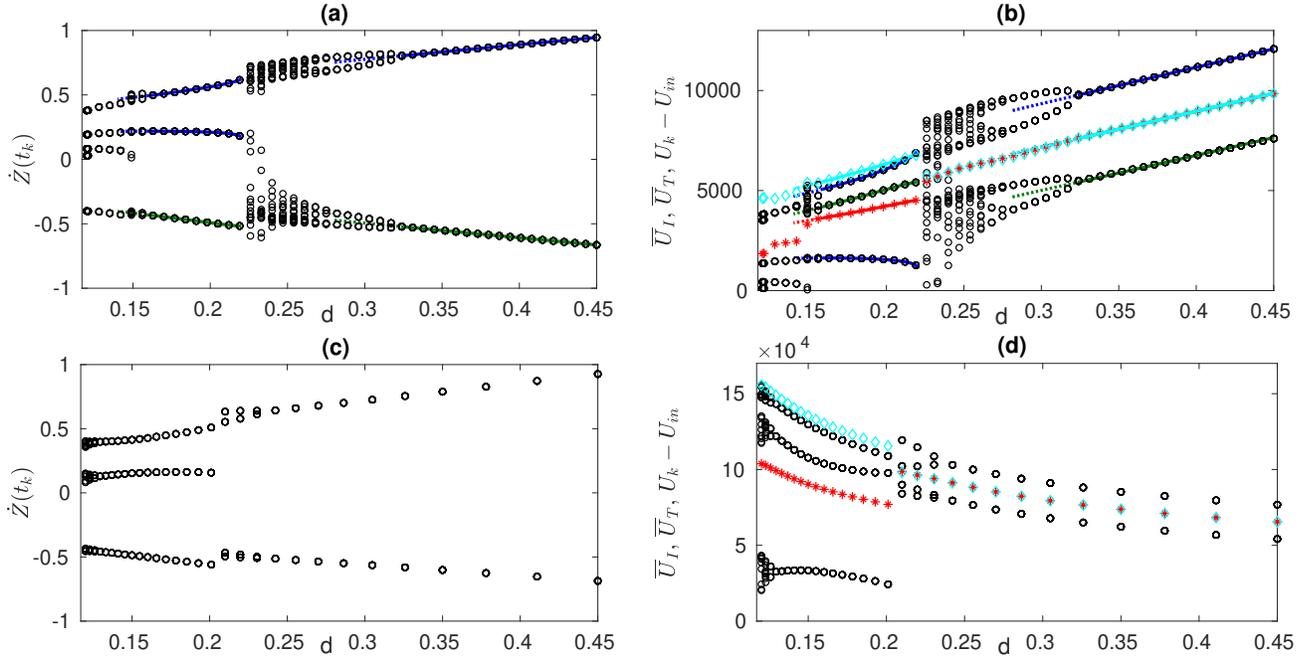}}
\caption{\scriptsize{Numerical (open circles $o$'s, stars $*$'s, diamonds $\lozenge$'s) and analytical stable/unstable (solid/dashed lines) values for impact velocities and output voltages for $\beta=\pi /3$. (a) Impact velocities (blue/green lines for bottom/top) for $0.19<s<0.72$, $\parallel \hat{F}\parallel = 5$.  The branches for the 2:1 solutions give, from top to bottom,  $\dot{Z}$ following the $P_3, P_1, P_2$ transitions.   (b) Output voltage $U_{k}-U_{\rm in}$ and average value of output voltage per impact $\overline {U}_I$ (red stars) and per unit time $\overline {U}_T$ (cyan diamonds) corresponding to $\dot Z$ in (a). The branches for the 2:1 solutions give, from top to bottom,  $U_k$ following the $P_3, P_2, P_1$ transitions; (c)-(d) Impact velocity and output voltage for $s=0.85$ with varying $\parallel \hat{F}\parallel$ between 6 and 22.}} \label{Fig0}
\end{figure}

In Figure \ref{Fig0} (a)-(d) we show the analytical and numerical results for the relative velocity at impact $\dot{Z}_k$ and output voltage  $U_k- U_{\rm in}$ vs. $d$ for the 1:1 periodic motion based on the results from \cite{Serdukova2019}. In the top row, the different values of  $d$ follow from variation in $s$, while in the bottom row the different values of $d$ follow from variation in $\hat{F}$. For decreasing values of $d$, the 1:1 period-2 motion loses stability via a sequence of period doubling bifurcations and eventually, for some parameter combinations, an apparently chaotic motion is observed for a window of values of $d$. For smaller values of $d$  these 1:1 behaviors with longer period or chaos are displaced by asymmetric motions with multiple impacts of $\partial B$ per period, that is, n:1 periodic motion. The analytical results corresponding to the branches for the 1:1 period-2 solutions are obtained both numerically and analytically in \cite{Serdukova2019}. Branches corresponding to period doubled 1:1 solutions and chaotic behavior are obtained numerically. The analytical results shown for the 2:1 period-2 motion are obtained in Sections 3 and 4 below, where we restrict our analysis to the derivation of the 2:1 solutions and their linear stability. Numerical results are compared to these analytical solutions, and also show additional bifurcations to n:1 periodic solutions, discussed further in Section 4. Analysis related to period doubling bifurcations, chaotic behavior, and grazing is left to future investigations.
\begin{figure}[b!]
 	\centering
 	\scalebox{0.59}
 	{\includegraphics{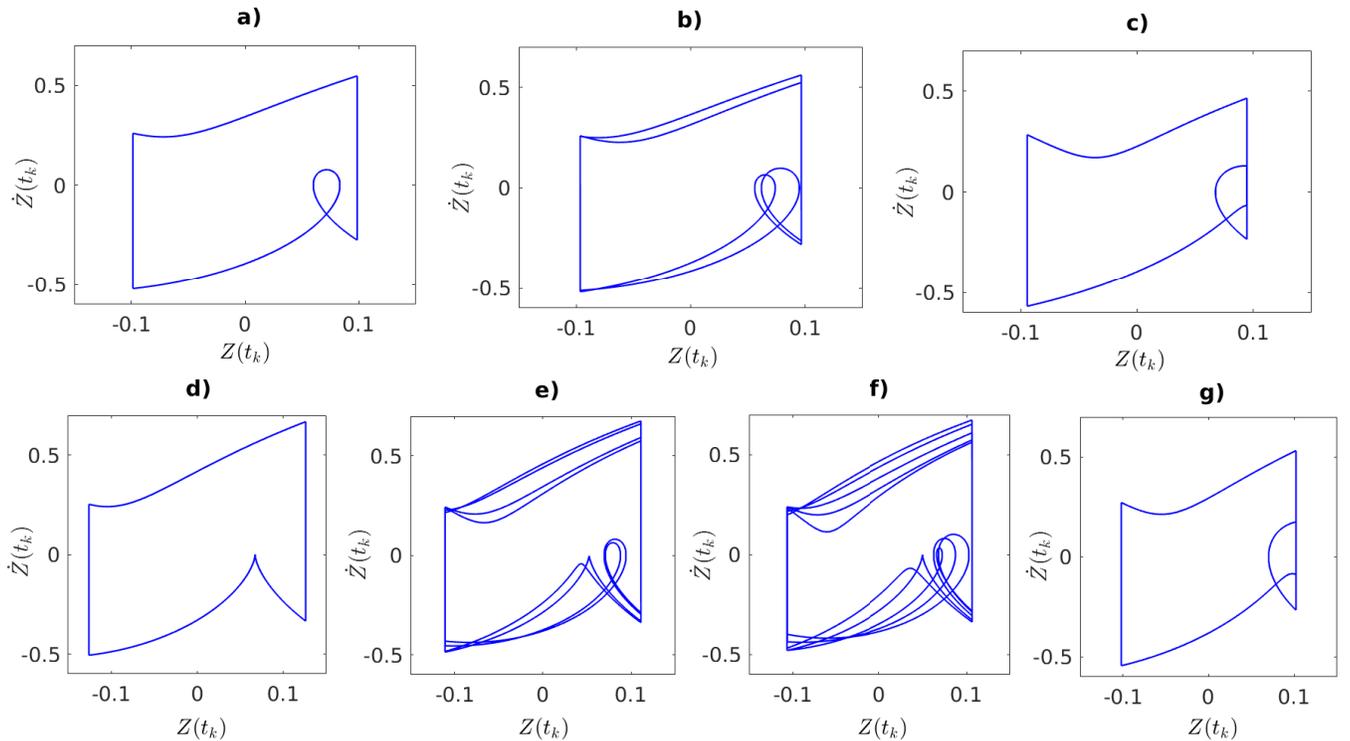}}
 	\caption{\scriptsize{Phase portrait for initial relative position $Z(t_0)= d/2$, $M = 124.5$ g, $r=0.5$. For (a)-(c) $\beta=\pi/2$,  $\parallel \hat{F}\parallel = 61$ N and $\omega = 18 \pi$ Hz; for (d)-(g) $\beta=\pi/6$,  $\parallel \hat{F}\parallel = 5$ N and $\omega = 5 \pi$ Hz. (a) 1:1 motion $d=0.197$, $s=0.316$,  $\dot{Z}(t_0)=0.5474$, $\varphi = 6.211$; (b) Grazing behavior for the 1:1 period-4 motion, with $d=0.193$, $s=0.309$,  $\dot{Z}(t_0)=0.561$, $\varphi = 6.229$; (c) 2:1 motion with $d=0.189$, $s=0.302$,  $\dot{Z}(t_0)=0.465$, $\varphi = 6.177$; (d) 1:1 motion $d=0.252$, $s=0.405$,  $\dot{Z}(t_0)=0.669$, $\varphi = 0.128$; (e) 1:1 period-8 motion $d=0.222$, $s=0.357$,  $\dot{Z}(t_0)=0.676$, $\varphi = 0.242$; (f) Grazing behavior of  1:1 period-10 motion $d=0.213$, $s=0.342$,  $\dot{Z}(t_0)=0.674$, $\varphi =0.321$; (g) 2:1 motion $d=0.204$, $s=0.328$,  $\dot{Z}(t_0)=0.532$, $\varphi =6.106$. } }\label{FigFN2}
 \end{figure} 
  
Figure \ref{Fig0} (b) and (d) shows the corresponding output voltage $U_{k}- U_{\rm in}$, for the same range of $d$ as in (a) and (c). Two different averaged output voltages are also shown, average per impact $\overline {U}_I$ and average per time unit $\overline {U}_T$, based on $30$ (non-dimensionalized) time units in $t$  ($\tau =6$ sec.) for 1:1 motion and 20 time units in $t$ ($\tau = 4$ sec.) for 2:1 motion. Note that the transition to different n:1 solutions corresponds to jumps in $\overline {U}_I$ and $\overline {U}_T$, given the change in the nature of the periodic solution. The additional impacts have low velocity, following naturally from the fact that they are born via grazing bifurcations, at which $\dot{Z}_j=0$ and $Z_j =  d/2$. For example, at $d=d_{\rm graz}$, $\dot{Z}_j=0$ and $Z_j = \pm d/2$, and there is a transition to 2:1 motion for $d<d_{\rm graz}$. Then the averaged per impact output voltage $\overline {U}_I$ drops for increasing $n$. The averaged output per unit time $\overline {U}_T$ is more complex, since the impact velocities following $P_{l}$ for $l=2, \; 3$ change with the addition of a low impact velocity from $P_1$. The increase in output voltage is achieved through increased cylinder length $s$, with other parameters fixed, or by increasing the forcing strength keeping $s$ constant, up to values of $d$ where there are transitions to n:1 periodic motions. The impact velocity and output voltage in Figure \ref{Fig0} (c) and (d) are obtained for  fixed $s=0.85$ and the variable strength of the forcing $6 < \parallel \hat{F}\parallel < 22$, in contrast to (a) and (b) for which $\parallel \hat{F}\parallel=5$ is fixed and $s$ varies for $0.118<d<0.448$. The nonlinear increase for  $U_{k}- U_{\rm in}$, $\overline {U}_I$ and $\overline {U}_T$ in (d) as $d$ decreases is due to the inverse dependency of  $d=\frac{s M \omega^2}{\parallel \hat{F} \parallel \pi^2}$ on the strength of the forcing. Note that the forcing $\bar{g}=\frac{M g \sin \beta}{\parallel \hat{F}\parallel}$  from the gravitational term  in \eqref{eq_relative} is also inversely proportional to $\parallel \hat{F} \parallel$.  
 
\begin{figure}[t!]
	\centering
	\scalebox{0.60}
	{\includegraphics {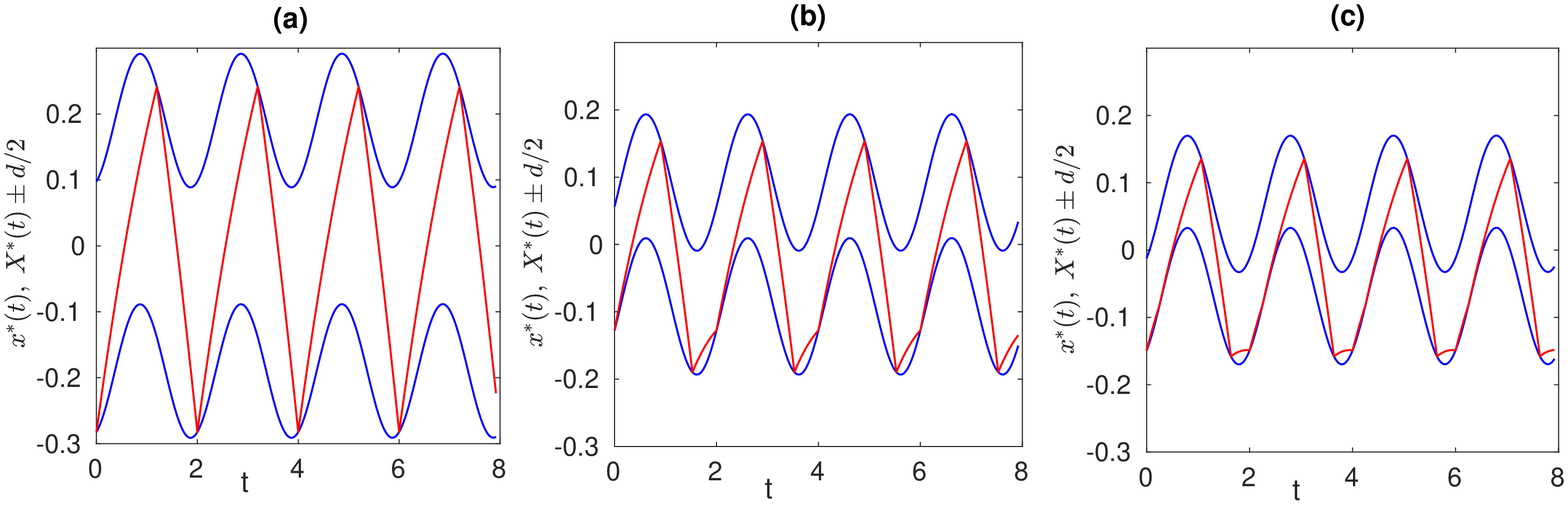}}
	\caption{\scriptsize{Time series of the period-2 absolute displacement of the  capsule top and bottom $X^*(t) \pm d/2$  (blue lines) and the absolute ball displacement  $x^*(t)$ (red line) for $t_0=0$ and $Z(t_0)= d/2$. (a) 1:1 motion for $d=0.38$, $s=0.61$, $\dot{Z}(t_0)=0.8673$, $\varphi = 0.4217$; (b) 2:1 motion for $d=0.184$, $s=0.30$,  $\dot{Z}(t_0)=0.2164$, $\varphi = 1.21$; (c) 3:1 motion for $d=0.137$, $s=0.22$,  $\dot{Z}(t_0)=0.2059$, $\varphi = 0.6503$. For all figures $M = 124.5$ g, $r=0.5$ and $\omega = 5 \pi$ Hz.}} \label{Fig0b}
\end{figure}

Figure \ref{FigFN2} illustrates the typical transition from 1:1 to 2:1 families of solutions in the phase plane, via a sequence of period doublings,  then grazing at a value of $d= d_{\rm graz}$ at which $\dot{Z}_k= 0$ and $Z_k = d/2$ as shown in Figure \ref{FigFN2} (b) and (f). From Figure \ref{Fig0} we see that 2:1 and other n:1 solutions persist for $d<d_{\rm graz}$, with $n$ increasing with decreasing $d$. To illustrate and compare the 1:1, 2:1 and 3:1 motions, Figure \ref{Fig0b} shows the absolute displacements $X^*(t)$ of the top and bottom of the cylinder under an external force and the motion of the ball $x^*(t)$ in the cylinder. These show the number of impacts per period in each case. 

 \section{Analytical expressions for periodic 2:1 motion}

In this section we obtain analytical expressions for the parametric dependence of the 2:1 period-2 motion, using the maps $P_1$, $P_2$ and $P_3$ for the sequence of impacts over the intervals $\Delta t_j$ for $j = k, k+1,k+2$. 

We derive equations for the quadruples $(\dot{Z}_{k},\varphi_{k},\Delta t_{k}, \Delta t_{k+1})$ corresponding to the 2:1 periodic solutions of \eqref{eq_relative} - \eqref{eq_impact}, in terms of the parameters $d$, $r$ and $\bar{g}$, with $\Delta t_{k}$ and $\varphi_{k} = {\rm mod}(\pi t_{k} + \varphi, 2\pi)$ as defined in Section \ref{sec:previous}.
We focus on a  2:1 period-T motion with three impacts per period $T$ of the forcing $f(t)$, so that
\begin{eqnarray}
t_{k+3}=T+t_k \; \;, Z_k= Z_{k+3} \; \;, \mbox{and} \; \; \dot{Z}_{k+3}=\dot{Z}_k \;. \label{period2}
\end{eqnarray}
The times for the transitions $P_1$, $ P_2$ and $P_3$ are defined as $T_1$, $T_2$ and $T_3$, with 
\begin{align}
\nonumber &T_1=\Delta t_{k}=t_{k+1}-t_k, \; \; T_2 = \Delta t_{k+1}=t_{k+2}-t_{k+1},\\ 
&T_3 = \Delta t_{k+2}=t_{k+3}-t_{k+2}, \; \; T = T_1+T_2+T_3. \label{T1T3}
\end{align}
The 2:1 period-T model is then described by the three maps $P_1$, $P_2$ and $P_3$ from \eqref{zdotj_p1} and \eqref{zj_p1}
\begin{align}
&P_1 \; : \; (Z_k\in \partial B, \dot{Z}_k, t_k) \; \mapsto \; (Z_{k+1}\in \partial B, \dot{Z}_{k+1}, t_{k+1}),\nonumber\\
&\dot{Z}_{k+1}=-r \dot{Z}_k + \bar{g}T_1 + F_1(t_{k+1})-F_1(t_k),\label{zdot_p1}\\
&0=-r \dot{Z}_k T_1 +\frac{\bar{g}}{2}T_1^2 +F_2(t_{k+1})-F_2(t_k)-F_1(t_k) T_1.\label{z_p1}
\end{align}
\begin{align}
&P_2 \; : \; (Z_{k+1}\in \partial B, \dot{Z}_{k+1}, t_{k+1}) \; \mapsto \; (Z_{k+2}\in \partial T, \dot{Z}_{k+2}, t_{k+2}),\nonumber\\
&\dot{Z}_{k+2}=-r \dot{Z}_{k+1}+\bar{g}T_2 +F_1(t_{k+2})-F_1(t_{k+1}),\label{zdot_p2} \\
&-d=-r \dot{Z}_{k+1}T_2 +\frac{\bar{g}}{2}T_2^2 +F_2(t_{k+2})-F_2(t_{k+1})-F_1(t_{k+1}) T_2. \label{z_p2}
\end{align}
\begin{align}
&P_3 \; : \; (Z_{k+2}\in \partial T, \dot{Z}_{k+2}, t_{k+2}) \; \mapsto \; (Z_{k+3}\in \partial B, \dot{Z}_{k+3}, t_{k+3}), \nonumber\\
&\dot{Z}_{k+3}=-r \dot{Z}_{k+2}+\bar{g}T_3 +F_1(t_{k+3})-F_1(t_{k+2}), \label{zdot_p3} \\
&d=-r \dot{Z}_{k+2}T_3 +\frac{\bar{g}}{2}T_3^2 +F_2(t_{k+3})-F_2(t_{k+2})-F_1(t_{k+2}) T_3. \label{z_p3}
\end{align}
We first use a number of substitutions to eliminate $\dot{Z}_{k+1}$, $\dot{Z}_{k+2}$ from \eqref{zdot_p1} - \eqref{z_p3} and obtain four equations in terms of $\dot{Z}_k$, from which we obtain $(\dot{Z}_{k},\varphi_{k},\Delta t_{k}, \Delta t_{k+1})$.
   
By adding \eqref{zdot_p1}, \eqref{zdot_p2}, \eqref{zdot_p3} and using the relationships $T = T_1+T_2+T_3$,  and $ F_1(t_{k+3})=F_1(T+t_k)=F_1(t_k)$, we obtain    
\begin{align}
\dot{Z}_k=\dfrac{1}{1-r+r^2} \left[(r-1)\bar{g}T_1 -\bar{g}T_2 + (1-r)F_1(t_k)+ r F_1(t_{k+1})-F_1(t_{k+2}) + \dfrac{T\bar{g}}{r+1}\right].\label{z0}
\end{align}
A second equation for $\dot{Z}_k$ is obtained from \eqref{z_p1} 
\begin{align}
\dot{Z}_k= \dfrac{1}{r T_1} \left[F_2(t_{k+1})-F_2(t_k) \right]+ \dfrac{1}{2 r} \left[\bar{g} T_1 - 2 F_1(t_k) \right]. \label{F1tk}   
\end{align} 
Substituting \eqref{zdot_p1} into \eqref{z_p2} yields a third expression for $\dot{Z}_k$  
\begin{align}
\dot{Z}_{k}=\dfrac{1}{r} \left[ \bar{g}T_1 + F_1(t_{k+1})-F_1(t_k) \right]-\dfrac{1}{r^2 T_2} \left[ d + F_2(t_{k+2})-F_2(t_{k+1}) \right] -\dfrac{1}{2 r^2} \left[ \bar{g} T_2 - 2 F_1(t_{k+1}) \right]. \label{F1tk1}
\end{align} 
Finally, adding \eqref{z_p1}, \eqref{z_p2}, \eqref{z_p3} and using relationship $F_2(t_{k+3})=F_2(T+t_k)=F_2(t_k)$ gives a fourth equation for $\dot{Z}_k$
\begin{align}
\dot{Z}_k&=\dfrac{1}{r^3 T_3 - r^2 T_2 + r T_1} \left[ \dfrac{\bar{g}}{2} (T_1^2+T_2^2+T_3^2)+F_1(t_k)(-r^2 T_3+rT_2-T_1)\right]+ \label{F1tk2} \\ 
\nonumber &+ \dfrac{1}{r^3 T_3 - r^2 T_2 + r T_1}\left[F_1(t_{k+1})(r^2 T_3-r T_2 +r T_3-T_2)+r^2 \bar{g} T_1 T_3-r \bar{g} T_1 T_2 \right] + \\
\nonumber &+ \dfrac{1}{r^3 T_3 - r^2 T_2 + r T_1} \left[ - r \bar{g} T_2 T_3 - (1+r) T_3 F_1(t_{k+2})\right].  
\end{align}
Then we solve \eqref{z0} - \eqref{F1tk2} to obtain $(\dot{Z}_{k},\varphi_{k},\Delta t_{k}, \Delta t_{k+1})$, using the Matlab function \textit{vpasolve}. A specific choice of $f(t)=\cos(\pi t+\varphi)$ for which
\begin{align}
&F_1(t)=\frac{1}{\pi} \sin(\pi t + \varphi) \; \; \mbox{and} \; \;  F_2(t)=-\frac{1}{\pi^2} \cos(\pi t + \varphi),
\end{align}
provides specifics for the equations for  $(\dot{Z}_k,\varphi_k,\Delta t_k, \Delta t_{k+1})$. It is convenient to write the time intervals between impacts in terms of the parameters $q$ and $p$, that capture the fractions of the period of forcing corresponding to each of the three impacts as follows,
\begin{align}
T_1=2nq, \; \; T_2=2np, \; \; T_3=2n(1-q-p),\; \; \mbox{and} \; \; T=2n\, . \label{T1T3n} 
\end{align}
We take $n=1$ for which the period is $T=2$. Without loss of generality, we take $t_k=0$, so $\varphi_k =$ mod$(\varphi, 2\pi)$. Then the four equations \eqref{z0} - \eqref{F1tk2} take the form of $\dot{Z}_k$ as functions of $\phi$, $q$, and $p$ 
\begin{align}
\dot{Z}_k&=\dfrac{1}{1-r+r^2} \left[2nq(r-1)\bar{g} -2np\bar{g} + \dfrac{1-r}{\pi}\sin(\pi t_k+\varphi) + \dfrac{r}{\pi} \sin(\pi [t_k+2nq]+\varphi)\right]+ \label{z0n}\\
&+ \dfrac{1}{1-r+r^2} \left[-\dfrac{1}{\pi}\sin(\pi [t_k+2nq+2np]+\varphi) + \dfrac{2n\bar{g}}{r+1}\right], \nonumber
\end{align}
\begin{align}
\dot{Z}_k = \dfrac{1}{\pi r} \left[ n \pi q \bar{g} -\sin(\pi t_k+\varphi) -\dfrac{1}{2n\pi q}\cos(\pi [t_k+2nq]+\varphi)+ \dfrac{1}{2n\pi q}\cos(\pi t_k+\varphi)\right], \label{F1tkn}
\end{align}
\begin{figure}[p]
\centering
\scalebox{0.58}
{\includegraphics{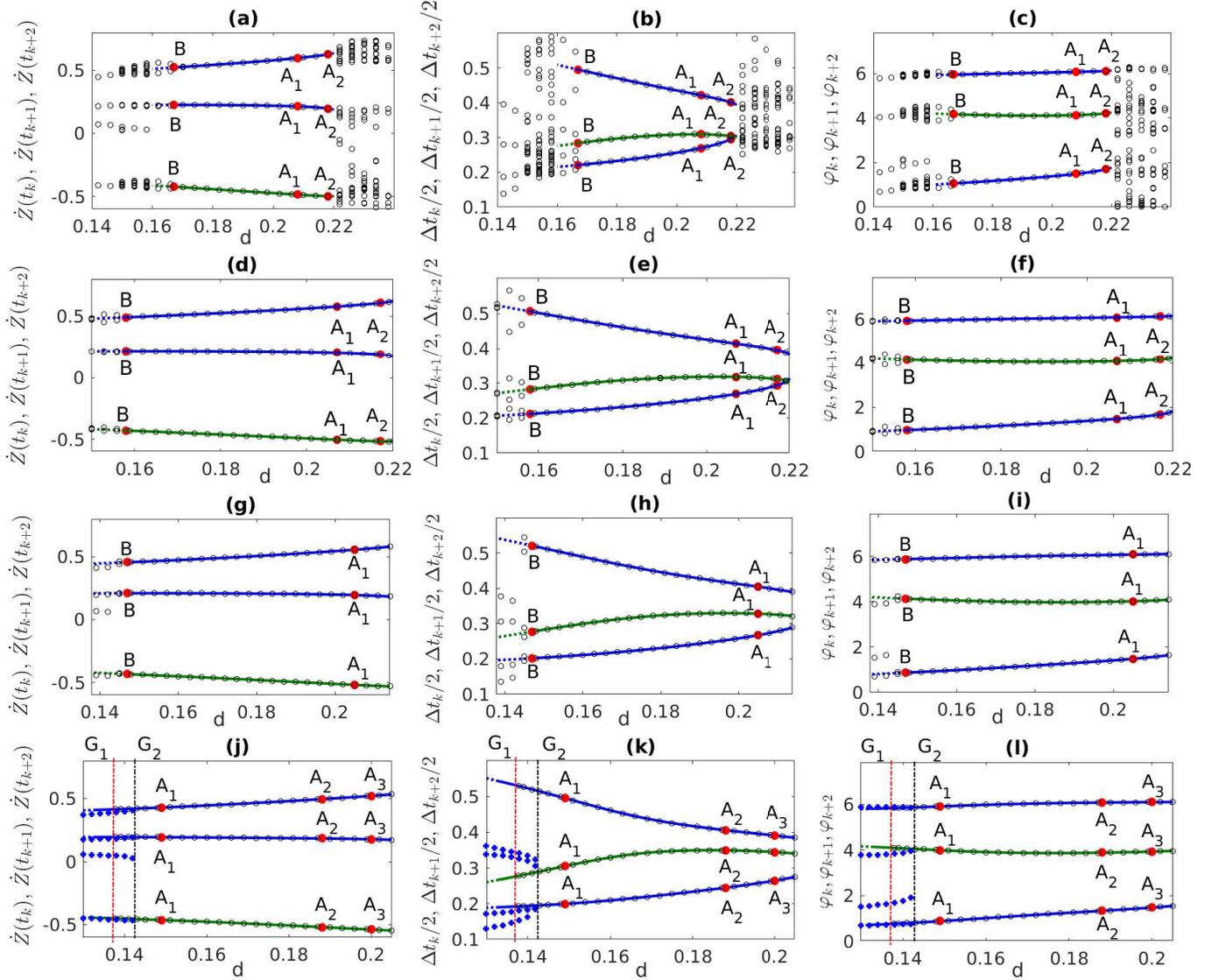}}
\caption{\scriptsize{Blue and green lines show the analytical results for the 2:1 periodic solutions, with numerical results indicated by open circles. Solid (dashed) lines correspond to stable (unstable) analytical solutions. (a)-(c) Asymmetric branches of the period-2 solutions for $\beta = \pi /2$, and $0.27 < s < 0.37$; (d)-(f) Asymmetric branches of the period-2 solutions for $\beta = \pi /3$ and $0.25 < s < 0.37$; (g)-(i) Asymmetric branches of the $2$-periodic solutions for $\beta = \pi /4$ and $0.22 < s < 0.33$; (j)-(l) Asymmetric branches of the period-2 solutions for $\beta = \pi /6$ and $0.22 < s < 0.33$. The vertical lines correspond to grazing bifurcations; $d=G_1$ (black) ($d=G_2$ (red)) for the transition from 3:1 to 2:1 (2:1 to 3:1) solutions with increasing (decreasing) $d$. In panels (a), (d), (g), (j) the branches for the 2:1 solutions give, from top to bottom,  $\dot{Z}_k$ at impacts following the $P_3, P_1, P_2$ transitions;  in panels (b), (e), (h), (k)  the branches for the 2:1 solutions give, from top to bottom,  $\Delta t_j/2$ for the $P_2, P_3, P_1$ transitions ; in panels (c), (f), (i), (l) the branches for the 2:1 solutions give, from top to bottom, the phase difference $\phi_k$ before the  $P_1, P_3, P_2$ transitions. For all figures $M = 124.5$ g, $r=0.5$, $\parallel \hat{F}\parallel = 5$ N and $\omega = 5 \pi$ Hz.}}
\label{Fig2}
\end{figure}
\begin{align}
\dot{Z}_k&= \dfrac{1}{\pi r^2} \left[ \sin(\pi [t_k+2nq]+\varphi) + 2n \pi q r\bar{g} + r \sin(\pi [t_k+2nq]+\varphi)- r \sin(\pi t_k+\varphi) \right] + \label{F1tk1n} \\
\nonumber &+ \dfrac{1}{\pi r^2} \left[ \dfrac{1}{2n\pi p}\cos(\pi [t_k+2nq+2np]+\varphi)-\dfrac{1}{2n\pi p} \cos(\pi [t_k+2nq]+\varphi)- n \pi p \bar{g} - \dfrac{\pi d}{2np} \right],
\end{align}
\begin{align} 
\dot{Z}_k &= \dfrac{\sin(\pi t_k+\varphi)(-2nr^2(1-p-q)+2npr-2nq)}{2nr^3(1-p-q) - 2npr^2 + 2nqr} - \label{F1tk2n} \\
\nonumber &- \dfrac{2 n \sin(\pi [t_k+2nq+2np]+\varphi)(1-p-q)(1+r)}{2n \pi r^3(1-p-q) - 2n\pi pr^2 + 2n\pi qr}\\ 
\nonumber &+ \dfrac{\sin(\pi [t_k+2nq]+\varphi)(2n r^2(1-p-q)-2npr +2nr(1-p-q)-2np)}{2nr^3(1-p-q) - 2npr^2 + 2nqr}+\\ 
\nonumber &+ \dfrac{4n^2r^2 \bar{g} q(1-p-q)-4n^2\bar{g}rpq - 4n^2\bar{g}rp(1-p-q)+\bar{g} (2n^2q^2+2n^2p^2+2n^2(1-p-q)^2)}{2nr^3(1-p-q) - 2npr^2 + 2nqr}.
\end{align}

Solving \eqref{z0n} - \eqref{F1tk2n} for varying $d$, one gets the quadruples $(\dot{Z}_{k}, \varphi_{k}, \Delta t_{k}, \Delta t_{k+1})$ for
2:1 period-2 solutions. Then $\dot{Z}_{k+1}$ is obtained from \eqref{zdot_p1} and substitution of \eqref{zdot_p1} into \eqref{zdot_p2} gives the equation for $\dot{Z}_{k+2}$ 
\begin{eqnarray}
\dot{Z}_{k+2}=r^2 \dot{Z}_k - r\bar{g}T_1 +\bar{g}T_2 + rF_1(t_k)- (1+r)F_1(t_{k+1}) +F_1(t_{k+2}). \label{z2}
\end{eqnarray}
 
Figure \ref{Fig2} shows the analytical solutions for these quadruples for different angles of incline $\beta$ and compares them to the values obtained from numerical simulations of equations \eqref{eq_relative} - \eqref{eq_impact}. The 2:1 period-2 solutions are stable only in the ranges of $0.167 <d< 0.22$ (a)-(c), $0.158 <d< 0.22$ (d)-(f), $0.147 <d< 0.214$ (g)-(i) and $0.1378 <d< 0.205$ (j)-(l). The stable 2:1 solutions, represented by the solid  blue lines (impacts on $\partial B$) and green lines (impacts on $\partial T$)  agree with the numerical solutions represented by black open circles.  The unstable 2:1 solutions represented by dashed lines are also found analytically. The points $A_1, A_2, A_3, B$  are the critical points that indicate a change in the type or stability or instability of the 2:1 solutions, based on the linear stability analysis. For the case of $\beta=\pi/6$ in the bottom row of Figure \ref{Fig2},  vertical lines  indicate the numerically detected grazing bifurcations at $d=G_1$ and $d=G_2$, corresponding to $Z_j=d/2$ and $\dot{Z}_j=0$. There are two different values, since the bifurcation value differs depending on whether it is obtained from decreasing the parameter $d$, yielding a transition from a 2:1 period-2 solution to a 3:1 period-2 solution at $d=G_1$, or by increasing $d$,  yielding a transition from 3:1 to 2:1 period-2 solutions at $d=G_2$. These results indicate a region of bi-stability for the 2:1 and 3:1 period-2 solutions, which we discuss briefly in Subsection \ref{grazing}. 

\section{Stability and Bifurcation of 2:1 period-2 motion}

\subsection{Linear stability analysis}

The critical points $A_j$, $B$ as shown in Figure \ref{Fig2} are obtained from a linear stability analysis around the quadruples $(\dot{Z}_{k},\varphi_{k},\Delta t_{k}, \Delta t_{k+1})$ corresponding to the asymmetric period-2 solutions. A complete review of this method can be found in \cite{Shaw1983,Luo2013,Luo2004}.

Considering a small perturbation $\delta \mathbf{H}_{k}$ to the fixed point $\mathbf{H}^*_{k} = (t_{k},\dot{Z}_{k})$, we obtain the equation for $\delta \mathbf{H}_{k+3}$ linearized about $\delta \mathbf{H}_{k}=0$,
\begin{eqnarray}
&\delta \mathbf{H}_{k+3} = \mathrm{D}P(\mathbf{H}^*_{k}) \delta \mathbf{H}_{k} = \mathrm{D}P_3(\mathbf{H}^*_{k+2}) \cdot \mathrm{D}P_2(\mathbf{H}^*_{k+1}) \cdot \mathrm{D}P_1(\mathbf{H}^*_{k}) \; \delta \mathbf{H}_{k}, \label{Lin1}
\end{eqnarray}
with
\begin{figure}[p]
\centering 
\scalebox{0.65}
{\includegraphics{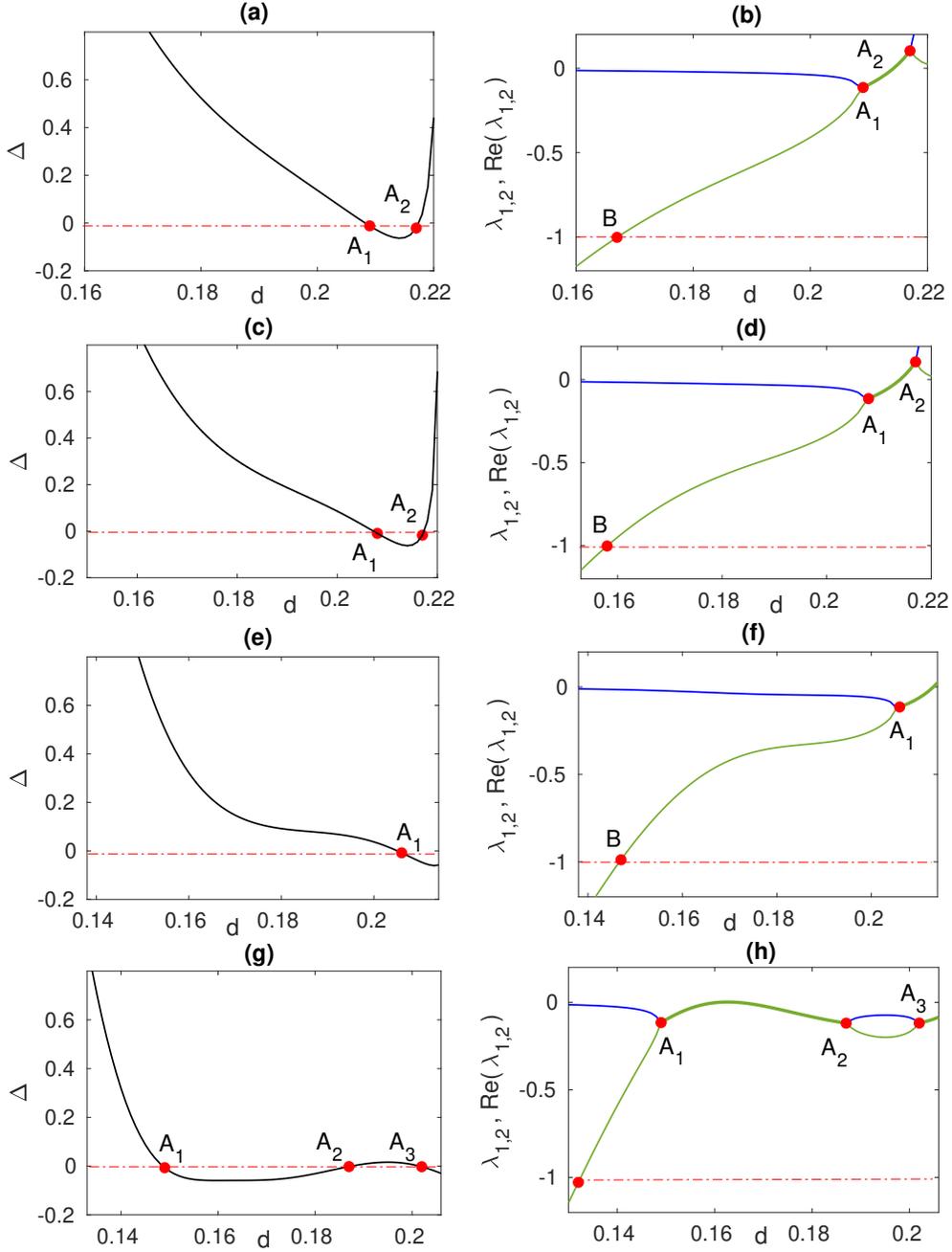}}
\caption{\scriptsize{Graphs of $\Delta$ (left column), and eigenvalues from the linear stability analysis (right column), showing real eigenvalues $\lambda_{1,2}$ (blue and green thin lines) and real part of complex eigenvalues ${\rm Re}\lambda_{1,2}$ (green thick line)  to confirm types and stability of solutions. In (a) and (b) for $\beta = \pi/2$, $0.27 < s < 0.37$; in (c) and (d) for $\beta = \pi /3$ and $0.25 < s < 0.37$; in (e) and (f) for $\beta = \pi /4$ and $0.22 < s < 0.33$; in (g) and (h) for $\beta = \pi /6$ and $0.22 < s < 0.33$. The red dot-dashed lines for $\Delta =0$ and $\lambda_{1,2}= - 1$ represent boundaries of the stability criteria. The left-most red circle in (h) corresponds to $\lambda_j=-1$ from the stability analysis. For all figures $M = 124.5$ g, $r=0.5$, $\parallel \hat{F}\parallel = 5$ N and $\omega = 5 \pi$ Hz.}}
\label{Fig3}
\end{figure}
\begin{align}
\mathrm{D}P = & \mathrm{D}P_3 \cdot \mathrm{D}P_2 \cdot \mathrm{D}P_1 = \nonumber \\ 
= &\left[
                                              \begin{array}{cc}
                                                      \frac{\partial t_{k+3}}{\partial t_{k+2}} & \frac{\partial t_{k+3}}{\partial \dot{Z}_{k+2}} \\
                                                   \frac{\partial \dot{Z}_{k+3}}{\partial t_{k+2}} & \frac{\partial \dot{Z}_{k+3}}{\partial \dot{Z}_{k+2}} \\
                                                  \end{array}
                                                 \right]_{\mathbf{H}_{k+2}=\mathbf{H}^*_{k+2}} 
                                                 \cdot & \left[
                                                   \begin{array}{cc}
                                                     \frac{\partial t_{k+2}}{\partial t_{k+1}} & \frac{\partial t_{k+2}}{\partial \dot{Z}_{k+1}} \\
                                                   \frac{\partial \dot{Z}_{k+2}}{\partial t_{k+1}} & \frac{\partial \dot{Z}_{k+2}}{\partial \dot{Z}_{k+1}} \\
                                                    \end{array}
                                                  \right]_{\mathbf{H}_{k+1}=\mathbf{H}^*_{k+1}} 
                                                 \cdot & \left[
                                                   \begin{array}{cc}
                                                     \frac{\partial t_{k+1}}{\partial t_{k}} & \frac{\partial t_{k+1}}{\partial \dot{Z}_{k}} \\
                                                   \frac{\partial \dot{Z}_{k+1}}{\partial t_{k}} & \frac{\partial \dot{Z}_{k+1}}{\partial \dot{Z}_{k}} \\
                                                    \end{array}
                                                  \right]_{\mathbf{H}_{k}=\mathbf{H}^*_{k}}.
                                                 \label{Lin2}
\end{align}
The entries $ \frac{\partial t_{l+1}}{\partial t_{l}} $, $ \frac{\partial t_{l+1}}{\partial \dot{Z}_{l}} $, $ \frac{\partial \dot{Z}_{l+1}}{\partial t_{l}} $, $ \frac{\partial \dot{Z}_{l+1}}{\partial \dot{Z}_{l}}$ for $l=k, k+1, k+2$ are given in A Appendix.

Using the trace ${\rm Tr}(DP)$ \eqref{Tr} and determinant ${\rm Det}(DP)$, the eigenvalues of the matrix $DP$ in (\ref{Lin2}) are computed by
\begin{align}
&\lambda_{1,2}=\frac{{\rm Tr}(DP)\pm \sqrt{\Delta}}{2}, \nonumber \\
&\Delta = [{\rm Tr}(DP)]^2-4{\rm Det}(DP), \label{Lambda}
\end{align}
and shown in Figure \ref{Fig3}. The corresponding stability and analytical bifurcation conditions as obtained from the linear stability analysis are described in Table \ref{tab} below.

\begin{table}[h!]
\centering
\begin{tabular}{l|l|l}
\hline\hline
\multicolumn {1}{c|}{\textbf{Interval}} & \multicolumn {1}{c|}{\textbf{Criteria}} & \multicolumn {1}{c}{\textbf{Stability}}\\[1ex]
\hline\hline
$d<d_{B}$ & $\Delta>0$ and $|\lambda_i|>1$ & unstable node\\
\hline
$d_{B}<d<d_{A_1}$,  $d_{A_2}<d<d_{A_3}$,  & $\Delta>0$ and $|\lambda_i|<1$ & stable node\\
\hline
$d_{A_1}<d<d_{A_2}$,  $d>d_{A_3}$ & $\Delta<0$ and $|\lambda_i|<1$ & stable focus\\
\hline\hline
\end{tabular}
\caption{Conditions for stability as obtained from the linear stability analysis and shown in Figure \ref{Fig0}, with, for example, $d_{A_j}$ corresponding to the value of $d$ at $A_j$.}
\label{tab}
\end{table}
In addition to these conditions, note that for smaller $\beta$, specifically $\beta= \pi/6$ in the last row of Figure \ref{Fig3}, the linear stability analysis indicates an eigenvalue  $\lambda_1<-1$ for $d<.133$. This  stability result is represented by a change from solid to dashed line for the analytical solutions shown in Figure \ref{Fig2} (j)-(l). However, the linear analysis does not capture the grazing bifurcations indicated by the dash-dotted vertical lines in Figure \ref{Fig2}. Then, in practice, the grazing bifurcation for $d>.133$, rather than the local linearized behavior, drives the transition from 2:1 to 3:1 period-2 solutions. The values of $d$ corresponding to grazing bifurcations are not included in Table \ref{tab}, but instead discussed in Subsection \ref{grazing} below.
\begin{figure}[t!]
\centering 
\scalebox{0.65}
{\includegraphics{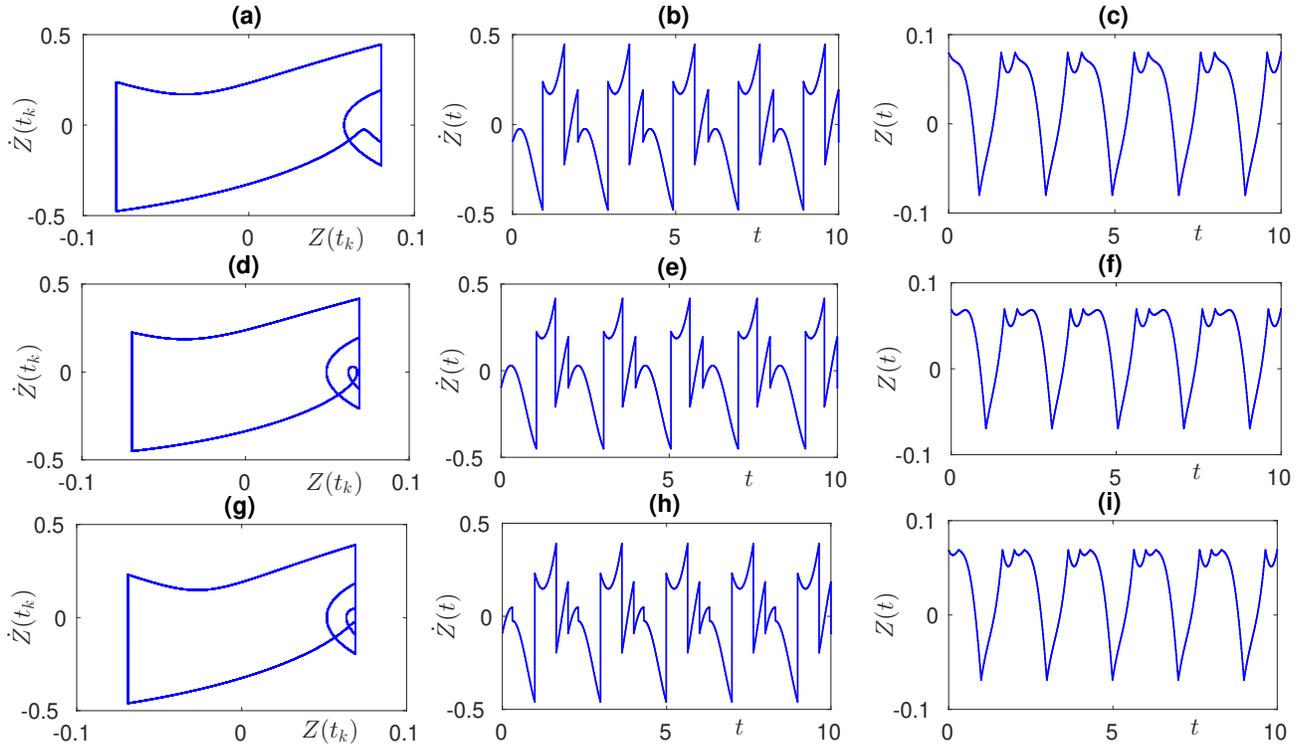}}
\caption{\scriptsize{Phase portrait and time series of period-2 motion, with $Z(t_0) = d/2$. (a)-(c) 2:1 motion for $\beta=\pi/6$, $d=0.16$, $s=0.27$,  $\dot{Z}(t_0)=0.1924$, $\varphi = 1.015$; (d)-(f) Grazing behavior of 2:1 motion for $\beta=\pi/6$, $d=0.139$, $s=0.23$, $\dot{Z}(t_0)=0.1959$, $\varphi = 0.7788$; (g)-(i) (3:1) motion for $\beta= \pi /6$, $d=0.138$, $s=0.23$, $\dot{Z}(t_0)=0.1845$, $\varphi = 0.7342$. For all figures $M = 124.5$ g, $r=0.5$, $\parallel \hat{F}\parallel = 5$ N and $\omega = 5 \pi$ Hz.}}
\label{Fig4}
\end{figure} 

If $\Delta<0$, as shown for $d_{A_1}< d <d_{A_2}$, $d>d_{A_3}$ and  in Figures \ref{Fig3} (a), (c), (e), (g), the eigenvalues of the matrix $DP$ are two complex conjugates. Their corresponding real parts $\rm Re(\lambda_i)= \rm Tr(DP)/2$ are shown in Figures \ref{Fig3} (b), (d), (f), (h), depicted by the thick green line. In these intervals the 2:1 period-2 solution is a stable focus since the eigenvalues also satisfy the condition $|\lambda_i|=\sqrt{\rm Det(DP)}<1$. 

If $\Delta>0$ and $\min_{i=1,2} (\lambda_i) < -1$, as in $d <d_{B}$ ranges in Figures \ref{Fig3} (b), (d), (f), the period-2 solution is an unstable node. The corresponding  critical point $B$ is a period doubling bifurcation. For the angles of incline $\beta = \pi/2$ and $\beta = \pi/3$ the stability behavior of the periodic motion is very similar revealing the predominance of node stability in the observed range of $d$ and having critical points of the same type: $B$ period doubling bifurcation, $A_1$ node/focus inflection and $A_2$ focus/node inflection. For smaller $\beta$, the qualitative behavior of the 2:1 period-2 solutions changes; specifically, grazing bifurcations drive the transition to 3:1 period-2 solutions for larger values of $d$ as compared with other critical values obtained from the linear stability analysis. We note that grazing bifurcations of the 2:1 period-2 solutions are observed for larger values of $\beta$ as well. They are not shown here since they occur for values of $d<d_B$ in those cases.

\subsection{The grazing transition and bistability}\label{grazing}

\begin{figure}[h!]
	\centering 
	\scalebox{0.56}
	{\includegraphics{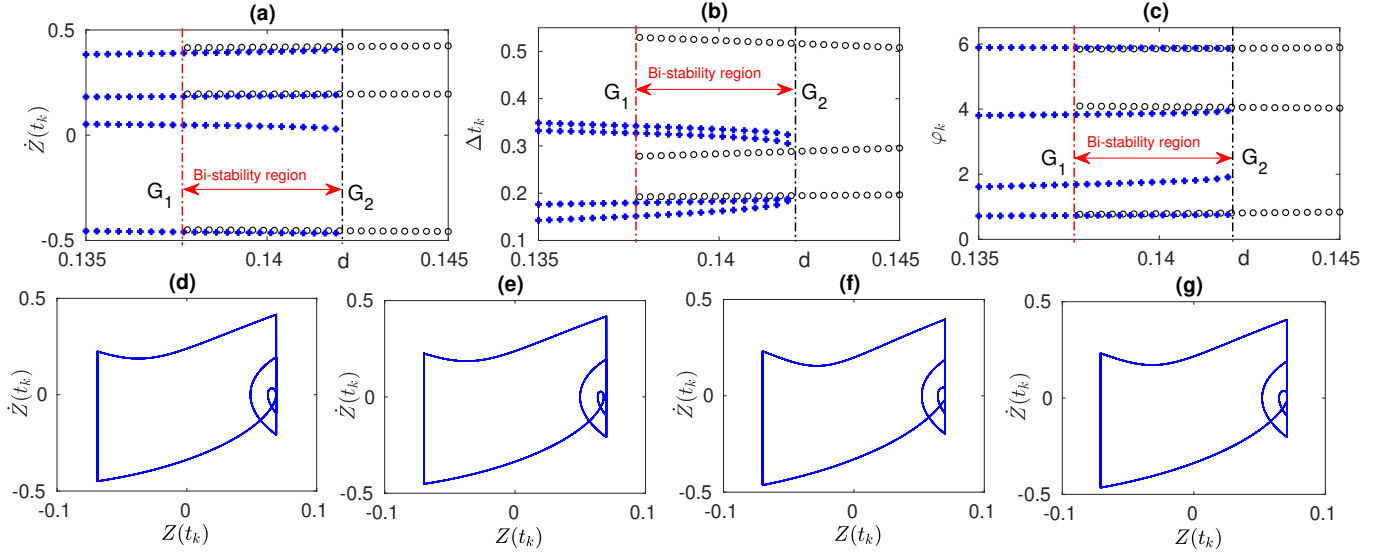}}
	\caption{\scriptsize{Bistable behavior of periodic 2:1 (black open circles o's) and 3:1 (blue crosses $+$'s) solutions in the vicinity of the grazing bifurcation $0.1378<d<0.1419$ for a) $\dot{Z}_{k}$, b) $\varphi_{k}$, c) $\Delta t_{k}$. Phase portraits with $Z(t_0)= d/2$ for d) Grazing point $G_1$ for $\beta=\pi/6$, $d=0.1378$, $s=0.221$,  $\dot{Z}(t_0)=0.416$, $\varphi = 5.842$; e) 2:1 motion for $\beta=\pi/6$, $d=0.14$, $s=0.224$,  $\dot{Z}(t_0)=0.4185$, $\varphi = 5.855$; f) 3:1 motion for $\beta=\pi/6$, $d=0.14$, $s=0.224$, $\dot{Z}(t_0)=0.3967$, $\varphi = 5.88$; g) Grazing point $G_2$ for $\beta=\pi/6$, $d=0.1419$, $s=0.228$,  $\dot{Z}(t_0)=0.4069$, $\varphi = 5.864$.}}
	\label{FigBiS}
\end{figure} 

For $\beta = \pi/6$ we numerically detect a different type of critical point for the 2:1 period-2 solutions, namely, grazing bifurcations as indicated by the vertical lines at $d=G_1$ and $d=G_2$ in Figure \ref{Fig2} (j)-(l), at which  $\dot{Z}_j = 0$ and $Z_j = d/2$ \cite{DiBernardo2008,Simpson2013,Simpson2018}. Figure \ref{FigBiS} zooms in on the bifurcation branches near these values. At these values of $d$ there are transitions between 2:1 and 3:1 period-2 motions. The transition from 2:1 to 3:1 period-2 behavior at $d=G_1$ is illustrated by the phase portrait and time series in Figure \ref{Fig4}. The initial conditions for these numerical simulations are obtained from the analytical expressions \eqref{z0n} - \eqref{F1tk2n}. In Figure \ref{Fig4} (d) the transition $P_2$ takes the form of a loop in the $\dot{Z}$ vs. $Z$ phase plane. As $d$ decreases, the loop intersects with $Z= d/2$, corresponding to an impact on $\partial  B$ with $\dot{Z}_{j} =0$. For decreasing $d$ this additional impact persists as shown in Figure \ref{Fig4} (d), yielding 3:1 period-2 solutions with an additional $P_1$ transition prior to $P_2$. 

Figure \ref{FigBiS} compares the grazing bifurcation at $d=G_1$ with a grazing bifurcation  that occurs as $d$ increases, leading to a transition from 3:1 to 2:1 period-2 solutions at $d=G_2$. The phase plane behavior for $d=G_1$ and $d=G_2$ are in panels (d) and (g), respectively. In addition, the bi-stability of 3:1 and 2:1 period-2 solutions for $G_1<d<G_2$ is shown via the bifurcation branches of $\dot{Z}_k$, $\varphi_k$ and $\Delta t_k$, as well as via the different phase plane behaviors at a $d=.14$ in this bistable region. While \cite{Luo2013} in chapter 6 explores some conditions for grazing and sticking and asymmetric behavior in the case with $\beta=0$, in general this bi-stability of different n:1 solutions via grazing has not been explored there or in other contexts. 

While not the focus of this paper, these results illustrate the importance of grazing bifurcations in driving different types of transitions in the VI-EH, as well as for the potential for hysteresis between bistable behaviors. The analytical conditions for this type of bifurcation in the case of the VI-EH is left for future investigation. 

\section{Energy output} 

Here we investigate the output voltage of the 2:1 period-2 behavior and compare these results with the 1:1 period-2 motion published in \cite{Serdukova2019}. Three variables corresponding to output voltage are shown, output voltage $U_k-U_{\rm in}$ at the $k^{\rm th}$ impact, average output per impact $\overline{U}_I$, and averaged output per unit of time $\overline{U}_T$. The derivation of $U_k-U_{\rm in}$ is summarized in \cite{Yurchenko2017} and $\overline{U}_I$, $\overline{U}_T$ are defined as
\begin{figure}[t!]
\centering 
\scalebox{0.59}
{\includegraphics{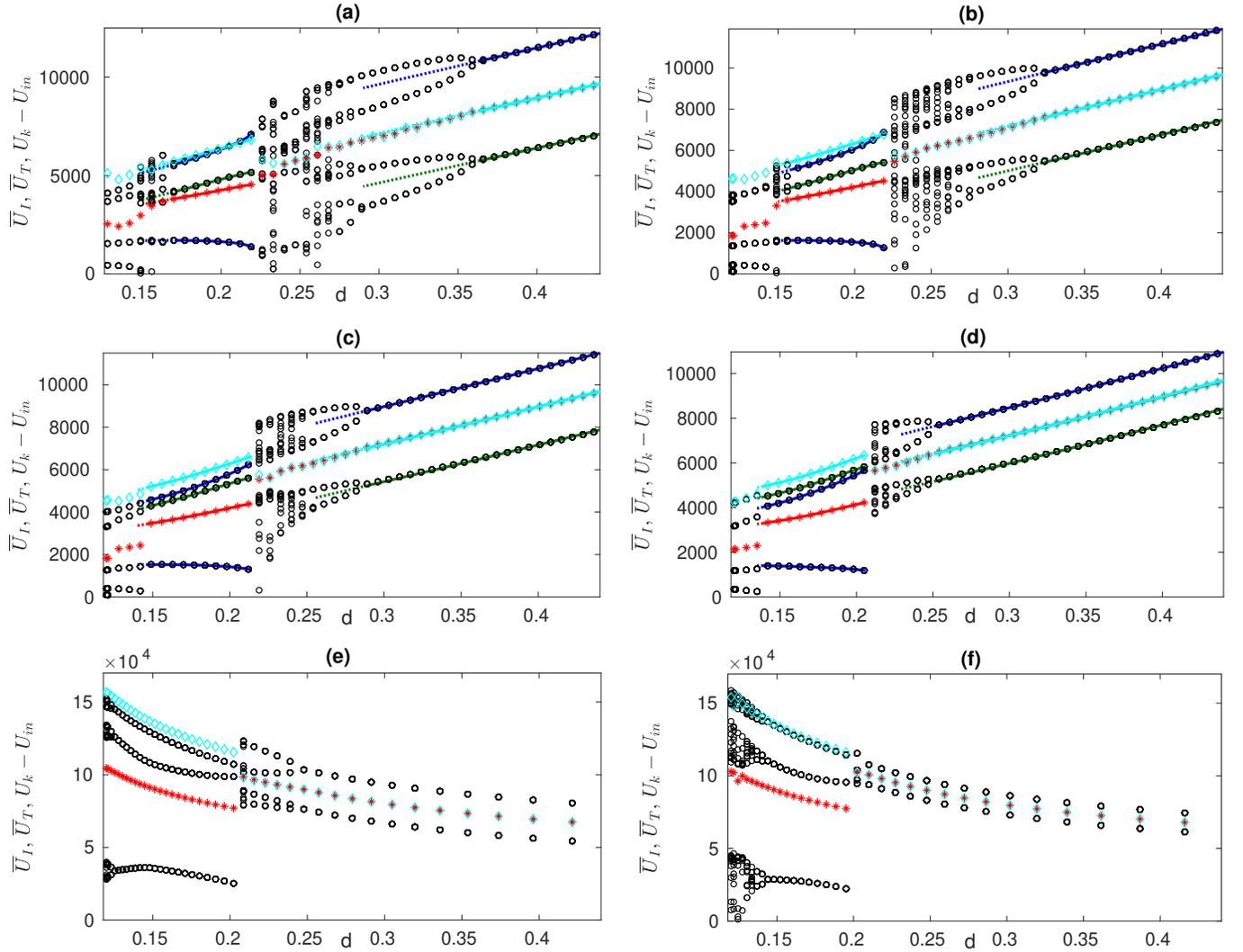}}
\caption{\scriptsize{Analytical results (solid and dashed lines) and numerical simulations (open circles $o$'s, stars $*$'s and diamonds $\lozenge$'s) for output voltage $\overline{U}_I$ (red)  and $ \overline{U}_T$ (cyan) and $U_k-U_{in}$ for (a) $\beta = \pi/2$, $0.19 < s < 0.72$; (b) $\beta = \pi /3$, $0.19 < s < 0.72$; (c) $\beta = \pi /4$, $0.19 < s < 0.72$; (d) $\beta = \pi /6$, $0.19 < s < 0.72$. For 2:1 period-2 solutions, in (a)-(c) the transitions $P_3, P_2, P_1$ are located from top to bottom, while in (d), $P_2, P_3, P_1$ are located from top to bottom. (e) For $\beta=\pi /2$, $s=0.85$ with varying $\parallel \hat{F}\parallel$ between 6 and 22. (f) For $\beta=\pi /6$, $s=0.85$ with varying $\parallel \hat{F}\parallel$ between 6 and 22. For all figures $M = 124.5$ g, $r=0.5$, $\omega = 5 \pi$ Hz.}}
\label{Fig5}
\end{figure}
\begin{eqnarray}
\overline{U}_I =\frac{\sum_{k=1}^N(U_k-U_{in})}{N}, \qquad
\overline{U}_T =\frac{\sum_{k=1}^N(U_k-U_{in})}{t_f - t_0},\label{V1}
\end{eqnarray}
where $N$ is the sample size of impacts and $t_f - t_0  = \frac{\omega}{\pi}(\tau_f-\tau_0) $ is the  corresponding non-dimensionalized time interval. We average over this time interval, since it is just a constant rescaling of the dimensionalized time interval, and then it is easy to compare $\overline{U}_I$ and $\overline{U}_T$ on the same plot. 

Figure \ref{Fig5} shows the output voltage for the 1:1 and 2:1 period-2 regimes, together with period doubled and chaotic regimes between these behaviors, for four different incline angles $\beta$. Panels (a)-(d) show variation due to  cylinder length  $s$ with fixed strength of forcing $\parallel \hat{F}\parallel$ and panels (e)-(f) show variations in $\parallel \hat{F}\parallel$ with fixed $s$. One obvious difference is the trend in output voltage, as observed previously in Figure \ref{Fig0}. Away from bifurcations, the output voltage increases with both increasing $\parallel F \parallel$ and increasing $s$. Then in (a)-(d) $U_k$ decreases with $d$ since $d$ is proportional to $s$, while in (e)-(f) $U_k$ shows a nonlinear increasing trend with decreasing $d$, due to the inverse relationship  $d=\frac{s M \omega^2}{\parallel \hat{F} \parallel \pi^2}$ to $\parallel \hat{F} \parallel$, as well as in the gravitational term $\bar{g}=\frac{M g \sin \beta}{\parallel \hat{F}\parallel}$.  
  
The bifurcations in the motion also result in changes in the output voltage, which we discuss in terms of the different measures of averaged output voltage. For 1:1 periodic motion, the average energy per impact $\overline{U}_I$ is equal to the average energy per unit of time $\overline{U}_T$, given that there are exactly two impacts for the 1:1  period-2 solutions. For the period doubled 1:1 solutions, as well as for more complex and chaotic behavior as shown for smaller values of $d>d_{\rm graz}$, we see a slight increase in the rate of decrease with $d$ of the average output voltage in (a)-(d), due primarily to the combination of values of impact velocities in the period doubled and more complex solutions that include some low velocity impacts. Following the  transition to 2:1 period-2 motion for $d<d_{\rm graz}$ the  average energy outputs $\overline{U}_I$ and $\overline{U}_T$ show jumps in the output value. Averaged output per impact $\overline{U}_I$ decreases due to the additional low velocity impact on $\partial B$ in the period $T=2$ for  2:1 period-2 solution. For the same reason, $\overline{U}_T$ increases due to this additional impact per period of the forcing. Similarly,  for the transition from 2:1 period-2 solutions to 3:1 period-2 solutions, the additional low velocity impact results in jumps both in $U_I$, which decreases across this critical value of $d$, and in $U_T$, which increases across this critical transition.  Note that here we show only the grazing transition at $d=G_1$ for 2:1 to 3:1 period-2 solutions, corresponding to decreasing $d$ in producing the bifurcation branches.
   
We also observe differences in the output voltages for different angles $\beta$ in terms of the  location in $d$ and sequence of period doubling bifurcations and complex or chaotic behavior, and for the value of $d$ at which the transition to 2:1 period-2 solutions occurs. In general, as $\beta$ increases, so do both the value of $d$ at which period doubling of the 1:1 solution occurs, and the value of $d_{\rm graz}$, the maximum value  for 2:1 period-2 solutions. Comparing Panels (a)-(d), for which $d$ decreases with $s$, and Panels (e)-(f), for which $d$ decreases with increasing $\parallel\hat{F}\parallel $, we observe a larger range of $d$ in (a)-(d) for period doubled and complex or chaotic behavior. This is partly due to the fact that even though $d$ decreases with increasing $\parallel \hat{F}\parallel $, the coefficient $\bar{g}$ also decreases with  increasing $\parallel \hat{F}\parallel $. Then for (e)-(f) as $d$ decreases there is a reduced influence of gravity, which would otherwise generate period doubled and complex behavior. For the 1:1 motion there is a small variation of the output voltage (less than 1\%) with $\beta$, for the maximum over the range of $d$ shown in panels (a)-(d), and similarly if we compare maximum output voltages over 2:1 motions for different angles. However, the parameter values at which these maxima occur differ with the incline $\beta$.
  
The result of this investigation suggests that the choice of the most efficient dynamical regime/device design in terms of the harvested electrical energy depends on the choice of measure for average output voltage and the changes in the parameter values of the system and the forcing.

\section{Conclusions}

In this paper we determine semi-analytical solutions and stability conditions for the 2:1 period-2 motion of an inclined vibro-impacting energy harvester (VI-EH). These results also provide insight into the VI-EH's energy harvesting potential. The device is  composed of a ball moving in a cylinder with dielectric elastomer (DE) material at the cylinder ends. It is driven by a harmonic forcing, and positioned with an incline angle. Energy is generated through impacts of the ball with the DE material, and the device exhibits n:m motion, where $n$ indicates  the number of impacts of the ball with the DE material on the bottom of the cylinder $\partial B$, and $m$ is the number of impacts on the top $\partial T$. Semi-analytical expressions for the generic period-T motion are derived through the three nonlinear mappings, that map the motion between the 3 impacts in the 2:1 motion per period.  These maps, together with  conditions that capture jump discontinuities in the velocity at impact, yield quadruples for the impact velocity, phase shift at impact, time intervals between the impacts. Analytical solutions are in excellent agreement with the numerical ones. Bifurcation points are obtained from a linear stability analysis around asymmetric periodic solutions. Based on the results it can be stated that: 

1. For larger values of the  incline angle $\beta$, the stability behavior of the 2:1 periodic motion exhibits predominance of node stability in the observed range of $d$. These solutions lose stability through  period doubling bifurcation for smaller values of $d$. This behavior is shown for $\beta = \pi/2$ and $\beta = \pi/3$.

2. For smaller values of incline $\beta$, the transition from 2:1 periodic behavior to 3:1 periodic behavior was observed as $d$ decreases. This transition occurs via a grazing bifurcation that is numerically detected. It occurs for larger values of $d$ compared with the values for other instabilities predicted by the linear analysis. These results are shown for $\beta= \pi/6$, for which bi-stability of the 2:1 and 3:1 solutions is numerically demonstrated near grazing.

3. The periodic asymmetric motions are less efficient compared to the motion with alternating top and bottom impacts per period of the forcing, when measured in terms of converted electrical energy per impact.

4. The 2:1 periodic motion results in  significant differences between the two measures of the harvested energy, averaged per impact, $U_{I}$, and averaged over time interval, $U_{T}$, giving greater value for $U_{T}$. Similar observations for 3:1 behavior are also shown.

\appendix 
\section{Appendix}

Here we give the details for the calculations of the eigenvalues $\lambda_{1,2}$. The entries in the matrices in \eqref{Lin2} are
\begin{align}
&\frac{\partial t_{k+1}}{\partial t_{k}} = \frac{r\dot{Z}_{k}-\bar{g}T_1-f(t_{k})T_1}{r \dot{Z}_{k}-\bar{g} T_1 - F_1(t_{k+1})+F_1(t_{k})}, \label{EntrP1}\\
&\nonumber\frac{\partial t_{k+1}}{\partial \dot{Z}_{k}} = \frac{-r T_1}{r \dot{Z}_{k}-\bar{g} T_1 - F_1(t_{k+1})+F_1(t_{k})},\\
&\nonumber \frac{\partial \dot{Z}_{k+1}}{\partial t_{k}} = \frac{\partial t_{k+1}}{\partial t_{k}}[f(t_{k+1})+\bar{g}]-[f(t_{k})+\bar{g}],\\
&\nonumber \frac{\partial \dot{Z}_{k+1}}{\partial \dot{Z}_{k}} = -r + \frac{\partial t_{k+1}}{\partial \dot{Z}_{k}} [f(t_{k+1})+\bar{g}],
\end{align}
\begin{align}
&\frac{\partial t_{k+2}}{\partial t_{k+1}} = \frac{r\dot{Z}_{k+1}-\bar{g}T_2-f(t_{k+1})T_2}{r \dot{Z}_{k+1}-\bar{g} T_2 - F_1(t_{k+2})+F_1(t_{k+1})}, \label{EntrP2}\\
&\nonumber\frac{\partial t_{k+2}}{\partial \dot{Z}_{k+1}} = \frac{-r T_2}{r \dot{Z}_{k+1}-\bar{g} T_2 - F_1(t_{k+2})+F_1(t_{k+1})},\\
&\nonumber \frac{\partial \dot{Z}_{k+2}}{\partial t_{k+1}} = \frac{\partial t_{k+2}}{\partial t_{k+1}}[f(t_{k+2})+\bar{g}]-[f(t_{k+1})+\bar{g}],\\
&\nonumber \frac{\partial \dot{Z}_{k+2}}{\partial \dot{Z}_{k+1}} = -r + \frac{\partial t_{k+2}}{\partial \dot{Z}_{k+1}} [f(t_{k+2})+\bar{g}],
\end{align}
and
\begin{align}
&\frac{\partial t_{k+3}}{\partial t_{k+2}} = \frac{r\dot{Z}_{k+2}-\bar{g}T_3-f(t_{k+2})T_3}{r \dot{Z}_{k+2}-\bar{g} T_3 - F_1(t_{k+3})+F_1(t_{k+2})}, \label{EntrP3}\\
&\nonumber\frac{\partial t_{k+3}}{\partial \dot{Z}_{k+2}} = \frac{-r T_3}{r \dot{Z}_{k+2}-\bar{g} T_3 - F_1(t_{k+3})+F_1(t_{k+2})},\\
&\nonumber \frac{\partial \dot{Z}_{k+3}}{\partial t_{k+2}} = \frac{\partial t_{k+3}}{\partial t_{k+2}}[f(t_{k+3})+\bar{g}]-[f(t_{k+2})+\bar{g}],\\
&\nonumber \frac{\partial \dot{Z}_{k+3}}{\partial \dot{Z}_{k+2}} = -r + \frac{\partial t_{k+3}}{\partial \dot{Z}_{k+2}} [f(t_{k+3})+\bar{g}].
\end{align}
For the period-2 motion the trace of the linearized matrix $DP$ are
\begin{align}
Tr(DP)&=-\dfrac{r^6 \dot{Z}(t_{k})}{F_1(t_{k+2})-F_1(t_{k+3})-r F_1(t_{k+1})+r F_1(t_{k+2}) +r^2 F_1(t_{k})-r^2 F_1(t_{k+1})+ \sigma_1} \label{Tr},
\end{align}
where $\sigma_1=r^3 \dot{Z}(t_{k})- \bar{g}T_3+r \bar{g} T_2-r^2 \bar{g}T_1$. The determinant of the linearized matrix $DP$ is a nonlinear function of $r$, $\bar{g}$, $\dot{Z}(t_{k})$, $T_1$, $T_2$, $T_3$, $f(t_{k})$, $f(t_{k+1})$, $f(t_{k+2})$, $f(t_{k+3})$, $F_1(t_{k})$, $F_1(t_{k+1})$, $F_1(t_{k+2})$ and $F_1(t_{k+3})$.

\bibliographystyle{ieeetr}
\bibliography{Bibliography1}
\end{document}